\input amstex
\documentstyle{amsppt}
\magnification=1200
\hcorrection{0.2in}
\vcorrection{-0.4in}

\define\volume{\operatorname{vol}}
\define\op#1{\operatorname{#1}}
\define\svolball#1#2{{\volume(\underline B_{#2}^{#1})}}
\define\svolann#1#2{{\volume(\underline A_{#2}^{#1})}}
\define\sball#1#2{{\underline B_{#2}^{#1}}}
\define\svolsp#1#2{{\volume(\partial \underline B_{#2}^{#1})}}

\NoRunningHeads
\NoBlackBoxes
\topmatter
\title Quantitative Maximal Diameter Rigidity of Positive Ricci
Curvature
\endtitle
\author Tianyin Ren \& Xiaochun Rong \footnote{Partially supported by NSFC 11821101, BNSF Z190003, and a research fund from Capital Normal University. \hfill{$\,$}\\ 2010 Mathematics Subject Classification. Primary 53C21, 53C23, 53C24.\hfill{$\,$}}
\endauthor
\address Mathematics Department, Rutgers University
New Brunswick, NJ 08903 USA
\endaddress
\email rong\@math.rutgers.edu
\endemail
\address Mathematics Department, Capital Normal University, Beijing,
P.R.C.
\endaddress
\email rtych12\@sina.com
\endemail
\abstract
In Riemannian geometry, the Cheng's maximal diameter rigidity theorem says that
if a complete $n$-manifold $M$ of Ricci curvature, $\op{Ric}_M\ge (n-1)$, has
the maximal diameter $\pi$, then $M$ is isometric to the unit sphere $S^n_1$.
The main result in this paper is a quantitative maximal diameter rigidity: if $M$
satisfies that $\op{Ric}_M\ge n-1$, $\op{diam}(M)\approx \pi$, and the
Riemannian universal cover of every metric ball in $M$ of a definite radius
satisfies a Riefenberg condition, then $M$ is diffeomorphic and bi-H\"older
close to $S^n_1$.
\endabstract
\endtopmatter
\document

\head 0. Introduction
\endhead

\vskip4mm

In Riemannian geometry, a maximal rigidity of Ricci curvature refers to a statement that a geometric or topological quantity of an $n$-manifold $M$ of $\op{Ric}_M\ge (n-1)H$ is bounded above by that of an $n$-manifold of constant sectional curvature $H$, and ``$=$'' implies that the sectional curvature on $M$ is constant $H$.

A quantitative maximal rigidity of Ricci curvature is a statement that if a geometric
quantity is almost maximal, then $M$ admits a nearby metric of constant sectional curvature
$H$. A significance of a quantitative maximal rigidity is that it implies the rigidity.

For a motivation of this paper, let's briefly recall three classical maximal rigidity results ($H=\pm 1, 0$) and their quantitative versions. A diffeomorphism between two Riemannian manifolds, $f: M_1\to M_0$, is called $\Psi(\epsilon|n)$-H\"older, (resp. bi-H\"older), if $f$ (resp. and $f^{-1}$) distorts distance, $\frac {d_0(f(x),f(y))}{d_1(x,y)^\alpha}\le e^{\Psi(\epsilon|n)}$, $0<\alpha<1$, where $\Psi(\epsilon |n)$ denotes a constant depending on $\epsilon$ such that $\Psi(\epsilon|n)\to 0$ as $\epsilon\to 0$ while $n$ is fixed.

\proclaim{Theorem 0.1} Let $M$ be a complete $n$-manifold of $\op{Ric}_M\ge (n-1)$.

\noindent {\rm (0.1.1)} {\rm (Maximal volume rigidity)} The volume, $\op{vol}(M)\le \op{vol}(S^n_1)$, and ``$=$'' if and only if $M$ is isometric to the unit $n$-sphere $S^n_1$.

\noindent {\rm (0.1.2)} {\rm (Quantitative maximal volume rigidity, \cite{CC1})} There  exists a constant $\epsilon(n)>0$ such that for $0<\epsilon\le \epsilon(n)$, if $\op{vol}(M)>\op{vol}(S^n_1)-\epsilon$, then $M$ is $\Psi(\epsilon|n)$ bi-H\"older diffeomorphic to $S^n_1$.
\endproclaim

A homeomorphism in (0.1.2) was obtained in \cite{Pe1}, that $M$ is close to $S^n_1$ in the Gromov-Hausdorff distance $d_{\op{GH}}$ was proved in \cite{Co1}, the optimal (0.1.2)  is proved in \cite{CC1} via the Reifenberg's method in geometric measure theory, and a canonical Reifenberg
method can be found in \cite{CJN}.

There is a local version of maximal volume rigidity for $H=\pm 1$ or $0$; the quantitative version for $H=1$ generalizes (0.1.2) (\cite{CRX1}).

The volume entropy of a compact Riemannian manifold $M$ is defined by
$$h(M)=\lim_{R\to \infty}\frac{\ln \op{vol}(B_R(\tilde p))}R,$$
where $\pi: (\tilde M,\tilde p)\to (M,p)$ is a Riemannian universal cover, $B_R(\tilde p)$
denotes the $R$-ball centered at $\tilde p$.

\proclaim{Theorem 0.2} Let $M$ be a compact $n$-manifold of $\op{Ric}_M\ge -(n-1)$.

\noindent {\rm (0.2.1)} {\rm (Maximal volume entropy rigidity, \cite{LW})}  The volume entropy $h(M)\le (n-1)$, and ``$=$'' implies that $M$ is isometric to a hyperbolic manifold.

\noindent {\rm (0.2.2)} {\rm (Quantitative maximal volume entropy rigidity, \cite{CRX1})} Given $n\ge 2, d>0$, there exists a constant, $\epsilon(n,d)>0$, such that for $0<\epsilon\le \epsilon(n,d)$, if $\op{diam}(M)\le d$ and $h(M)>n-1-\epsilon$, then $M$ is $\Psi(\epsilon|n,d)$ bi-H\"older diffeomorphic to a hyperbolic manifold.
\endproclaim

Note that (0.2.2) was stated in \cite{CRX1} that $M$ is diffeomorphic and Gromov-Hausdorff close to a hyperbolic manifold. However, from its proof (that $\tilde M$ is GH-close to $\Bbb H^n$ and $\op{vol}(M)\ge v(n,d)>0$) the bi-H\"older closeness is obvious.

Note that different from (0.1.2), (0.2.2) will be false without the additional condition, $\op{diam}(M)\le d$ (\cite{GT}). Moreover, (0.2.2) does not hold if one replaces ``$\op{diam}(M)\le d$'' by ``$\op{vol}(M)\le  v$'' (\cite{SS}).

\proclaim{Theorem 0.3}  \noindent {\rm (0.3.1)} {\rm (Maximal first Betti number rigidity, \cite{Bo})} If $M$ is a compact $n$-manifold of $\op{Ric}_M\ge 0$, then the first Betti
number, $b_1(M)\le n$, and ``$=$'' implies that $M$ is isometric to a flat torus.

\noindent {\rm (0.3.2)} {\rm (Quantitative maximal first Betti number rigidity, \cite{Co2}, \cite{CC2})} Given $n\ge 2$, there exists a constant, $\epsilon(n)>0$, such that for $0<\epsilon\le \epsilon(n)$, if $M$ of $b_1(M)=n$ satisfies that $\op{Ric}_M\ge -\epsilon$
and $\op{diam}(M)\le 1$, then $M$ is $\Psi(\epsilon|n)$ bi-H\"older diffeomorphic to a flat torus.
\endproclaim

Similar to (0.2.2), (0.3.2) does not hold without the restriction on diameter.  For further work generalizing (0.3.2) via smoothing method (which yields a bi-H\"older closeness in (0.3.2)), see \cite{HW2}, \cite{HKRX} and \cite{Ro3}.

A partial motivation for this paper is from the lack of a quantitative version to
the classical Cheng's maximal diameter rigidity below.

\proclaim{Theorem 0.4} Let $M$ be a complete
$n$-manifold of $\op{Ric}_M\ge (n-1)$.

\noindent {\rm (0.4.1)} {\rm (\cite{My})} $\op{diam}(M)\le \pi$, thus $\pi_1(M)$ is finite.

\noindent {\rm (0.4.2)} {\rm (Maximal diameter rigidity, \cite{Che})} If $\op{diam}(M)=\pi$, then $M$ is isometric to $S^n_1$.
\endproclaim

A possible reason is that given any $\epsilon>0$, $\Bbb CP^n$ ($n\ge 2$, \cite{An1}) and $S^k\times S^{n-k}$ ($k\ge 2$ and $n-k\ge 3$, \cite{Ot}) admits a metric $g_\epsilon$ such that $\op{Ric}_{g_\epsilon}\ge (n-1)$ and $\op{diam}(g_\epsilon)\ge \pi-\epsilon$. Hence, similar to (0.2.2)  or (0.2.3), additional conditions are required for a quantitative maximal diameter rigidity.

Instead, there have been topological stability results asserting $M$ a sphere, homeomorphic or diffeomorphic, under various additional restrictions; either $\op{sec}_M\ge -\kappa$  (\cite{Pe2}), or $\op{Ric}_M\le C$ and a local regularity that implies that $\op{vol}(M)$ is not small (\cite{GP}, \cite{Pet}, \cite{PSZ}, \cite{PZ}, \cite{Wu}, etc);
none of them implies (0.4.2). Here we present an optimal diameter sphere theorem of Perel'man.

\proclaim{Theorem 0.5} {\rm (Diameter spheres, \cite{Pe2})} Given $n\ge 2, \kappa>0$, there is a constant $\epsilon(n,\kappa)>0$, such that if a complete manifold $M$ satisfies
$$\op{Ric}_M\ge n-1, \quad \op{diam}(M)\ge \pi-\epsilon(n,\kappa),\quad \op{sec}_M\ge  -\kappa,$$ then $M$ is homeomorphic to a sphere.
\endproclaim

Theorem 0.5 is a topological stability result, and the geometry on $M$ may be far from that of $S^n_1$; for instance, $M$ in Theorem 0.5 may have a very small volume (e.g., $M$ is obtained via smoothing around the two vertices of a spherical suspension over a round $(n-1)$-sphere of radius $\epsilon$.).

The main purpose in this paper is to establish a quantitative maximal diameter rigidity (see Theorems A and B).

First, according to \cite{CC2} a necessary and sufficient local regularity for a sequence of $n$-manifolds, $M_i$ of $\op{Ric}_{M_i}\ge -(n-1)$, converging in the Gromov-Hausdorff
topology to a compact $n$-manifold $M$, $M_i@>\op{GH}>>M$, is that there are constants, $\rho=\rho(M),  \delta(n)>0$, such that for $i$ large, any $x_i\in M_i$ is $(\rho,\delta(n))$-Reifenberg point (briefly, RP) i.e.,
$$d_{GH}(B_s(x_i),\b B^n_s(0))< s\cdot \delta(n),\quad \forall\,\, 0<s\le \rho,$$
where $\b B^n_s(0)$ denotes a $s$-ball in $\Bbb R^n$.

In our quantitative maximal diameter rigidity theorem, the local regularity condition is weak: a point $x\in M$ is called a $(\rho,\delta)$-local rewinding Reifenberg point (briefly, LRRP), if $\tilde x$ (not $x$) is a $(\rho,\delta)$-Reifenberg point, where $\pi: (\widetilde {B_\rho(x)},\tilde x)\to (B_\rho(x),x)$ is the Riemannian universal cover map. Note that $x$ being
a $(\rho,\delta)$-LRRP is a local regularity of `a collapsed metric', because the volume ratio, $\frac{\op{vol}(B_\rho(x))}{\op{vol}(\underline B^n_\rho(0))}$, can be arbitrarily small.

Indeed, the notion of LRRP was introduced by the second author in his proposal to extend the nilpotent structures on collapsed manifolds with bounded sectional curvature, constructed by Cheeger-Fukaya-Gromov (\cite{CFG}), to collapsed manifolds with Ricci curvature bounded below and local bounded covering geometry. For recent progress, see \cite{CRX1,2}, \cite{DWY}, \cite{HKRX}, \cite{HRW}, \cite{HW1,2}, \cite{NT1,2}, \cite{NZ},  \cite{PWY}, \cite{Ro1-3}, etc.

We now begin to state the main results in this paper.

\vskip2mm

\proclaim{Theorem A} {\rm (Quantitative maximal diameter rigidity)} Given $n\ge 2,  1\ge \rho>0$, there exist constants, $\delta(n), \epsilon(n,\rho)>0$, such that if a complete $n$-manifold $M$ satisfies
$$\op{Ric}_M\ge (n-1),\quad\op{diam}(M)\ge \pi-\epsilon(n,\rho),\quad \forall\, x\in M \text{ is a $(\rho,\delta(n))$-LRRP},$$
then $M$ is $\Psi(\delta(n))$ bi-H\"older diffeomorphic to $S^n_1$. Moreover, for $0<\delta<\delta(n)$,
there exists $0<\Psi(\delta|n,\rho)<\epsilon(n,\rho)$, such that if all $x\in M$ are $(\rho,\delta)$-LRRP and $\op{diam}(M)\ge \pi-\Psi(\delta|n,\rho)$, then $M$ is $\Phi(\delta|n,\rho)$ bi-H\"older diffeomorphic to $S^n_1$.
\endproclaim

Theorem A easily implies (0.4.2) (see Section 1). Theorem A is false if one replaces ``$(\rho,\delta(n))$-LRRP'' with ``$\op{vol}(M)\ge v>0$'' (\cite{An1}, \cite{Ot}). By Theorem A, the local regularity (designed for collapsing), ``$(\rho,\delta(n))$-LRRP'', implies the (possibly sharp) local regularity, ``$(c(n,\rho),\delta(n))$-RP'' (see Corollary 1.3), based on the key volume estimate from below on $\op{vol}(B_\rho(x))$ (see Theorem 1.1), where $c(n,\rho)>0$ is
a constant depending on $n$ and $\rho$.

According to \cite{CRX1} (seen from the proof of Theorem A in \cite{CRX1}, where a
smoothing technique via Ricci flows is applied), if the Riemannian universal cover, $\tilde M$, of any complete $n$-manifold $M$ with $\op{Ric}_M\ge (n-1)$, is $\op{GH}$-close to $S^n_1$, then $M$ is bi-H\"older diffeomorphic to a spherical space form (here $\pi_1(M)$ may have any large order).

In view of the above, we state the following strong version of Theorem A.

\proclaim{Theorem B} {\rm (Strong quantitative maximal diameter rigidity)} Given $n\ge 2, 1\ge \rho>0$, there exist constants, $\delta(n), \epsilon(n,\rho)>0$, such that if the Riemannian universal covering $\tilde M$ of a complete $n$-manifold $M$ satisfies
$$\op{Ric}_{\tilde M}\ge (n-1),\quad\op{diam}(\tilde M)\ge \pi-\epsilon(n,\rho),\quad \forall\, \tilde x\in \tilde M \text{ is a $(\rho,\delta(n))$-LRRP},$$
then $M$ is $\Psi(\delta(n))$ bi-H\"older diffeomorphic to a spherical space form. Moreover, for $0<\delta<\delta(n)$,
there exists $0<\Psi(\delta|n)<\epsilon(n,\rho)$, such that if all $\tilde x\in \tilde M$ are $(\rho,\delta)$-LRRP and $\op{diam}(\tilde M)\ge \pi-\Psi(\delta|n)$, then $M$ is $\Phi(\delta |n,\rho)$ bi-H\"older diffeomorphic to spherical space form.
\endproclaim

Note that every spherical $n$-space form satisfies Theorem B, among them only $S^n_1$ satisfies Theorem A.


By the non-collapsing property in Theorem 1.1, one is able to conclude two versions of quantitative maximal diameter rigidity results, which
assert that if local universal cover of $\rho$-balls on $M$ satisfies a regularity stronger than LRRP, then $M$ is diffeomorphic to a spherical space form with a corresponding strong global regularity.

In the case that $C\ge \op{Ric}_{\tilde M}\ge n-1$, we have the following
quantitative maximal diameter rigidity with $C^{1,\alpha}$-regularity (cf. \cite{An2}).

\proclaim{Corollary 0.6} {\rm(Diffeomorphic and $C^{1,\alpha}$ close to spherical space forms)} Let $M$ be as in Theorem B. If, in addition, $\op{Ric}_{\tilde M}\le C$, then $M$ is $\Psi(\epsilon |n,C)$ $C^{1,\alpha}$ diffeomorphic to a spherical space form ($0<\alpha<1$).
\endproclaim

Note that different from Theorem B, $\Psi(\epsilon|n,\rho)$ is independent of $\delta$; because for any $\delta<\delta(n)$, $\tilde X$ is a $C^{1,\alpha}$-manifold, and $\op{diam}(\tilde M_i)\to \pi$ ($\epsilon\to 0$) implies a suspension structure on $\tilde X$, thus $\tilde X$ is isometric to $S^n_1$ (see the proof of Theorem A by assuming Theorem 1.1).

For $x\in M$, we call the number, $\op{injrad}_\rho(\tilde x):=\min\{\op{injrad}(\tilde x),\rho\}$, the $\rho$-rewinding injectivity radius, where $\pi: (\widetilde {B_\rho(x)},\tilde x)\to (B_\rho(x),x)$ is
the universal covering map.

A quantitative maximal diameter rigidity with a regularity between Theorem A and Corollary 0.6
is the following (cf. \cite{AC}):

\proclaim{Corollary 0.7} {\rm(Diffeomorphic and $C^\alpha$-close to spherical space forms)} Given $n\ge 2, \rho, i_0>0$, there exist constants, $\delta(n), \epsilon(n,\rho, i_0)>0$, such that for $0<\epsilon\le \epsilon (n,\rho,i_0)$, if the Riemannian universal covering $\tilde M$ of a complete $n$-manifold $M$ satisfies
$$\op{Ric}_{\tilde M}\ge (n-1),\quad \op{diam}(\tilde M)\ge \pi-\epsilon, \quad \op{injrad}_\rho(\tilde x)\ge i_0, \quad \forall\, \tilde x\in \tilde M,$$
then $M$ is $\Psi(\epsilon |n,\rho, i_0)$ $C^\alpha$ diffeomorphic
to a spherical space form.
\endproclaim

Consider a compact $n$-manifold $M$ of $\op{Ric}_M\ge (n-1)$ and $\op{diam}(M)>\pi-\epsilon$, $p, q\in M$, $d(p,q)=\op{diam}(M)$. A consequence of the Ricci curvature-diameter condition is that the excess function,
$e_{p,q}(x)=d(p,x)+d(x,q)-d(p,q)$, $x\in M$, is small (\cite{AG}, \cite{GP}). Clearly, one may relaxes the Ricci curvature-diameter condition and requires a small excess, one may expect $M$ a homeomorphic sphere, provided additional restrictions: $\pi_1(M)$ is finite, a lower bound either on sectional curvature, or on the injectivity radius (\cite{GP}, \cite{PZ}, \cite{Wu}, etc.).

Our proof of Theorem A easily extends to the following excess sphere theorem.

\proclaim{Theorem C} {\rm (Excess spheres)} Given $n\ge 2, \rho, d>0$, there exist constants,
$\delta=\delta(n), \epsilon=\epsilon(n,\rho, d)>0$, such that if a compact $n$-manifold $M$ of a
finite fundamental group satisfies
$$\op{Ric}_M\ge -(n-1), \,  e_{p,q}(x)<\epsilon,\, d(p,q)=\op{diam}(M)\le d, \, \forall\, x\in  M \text{ is a $(\rho,\delta(n))$-LRRP},$$
then for $n\ge 4$, $M$ is homeomorphic to a sphere and $\op{vol}(M)\ge v(n,\rho,d)>0$, and for $n=3$, $M$ is homeomorphic to a spherical space form.
\endproclaim

The non-collapsing property in Theorem C for $n\ge 4$ is new, which has applications (see Remark 0.8), and Theorem C fails for $n=3$ (e.g., a sequence of lens spaces, $S^3_1/\Bbb Z_h, h\to \infty$, collapses to a round sphere, $S^2_{\frac 12}$, of radius $\frac 12$).

\remark{Remark \rm 0.8} In Theorem C (similar to Theorem A),  the weak local regularity, $(\rho,\delta(n))$-LRRP,  implies an `optimal'' local regularity $(c_1(n,\rho,d),\delta(n))$-RP ($n\ge 4$), thus Theorem C can have corollaries similar to Corollaries 0.6 and 0.7; here we omit precise statements.
\endremark

\vskip2mm

The rest of the paper is organized as follows:

\vskip2mm

In Section 1, we will prove Theorem A by assuming that $M$ is not collapsed (Theorem 1.1).  And we will outline
our approach to Theorem 1.1.

In Section 2, we will prove Theorem 1.1 modulo two technique lemmas, Lemmas 2.2 and 2.5.

In Section 3, we will prove Lemma 2.2.

In Section 4, we will prove Lemma 2.5.

In Section 5, we will prove Theorem C.

\vskip4mm

\head 1.  Proof of  Theorem A  (Non-collapsing Case)
\endhead

\vskip4mm

Our approach to Theorem A is to show that $M$ has an almost maximal volume, thus by (0.1.2) we conclude the desired result. In our proof, the main technical result is the following:

\proclaim{Theorem 1.1} {\rm(Volume estimate)} Given $n\ge 2, 1\ge \rho>0$, there exist constants, $\delta(n)$, $\epsilon(n,\rho)$, $v(n,\rho)>0$, such that if a complete $n$-manifold $M$ satisfies
$$\op{Ric}_M\ge (n-1),\quad\op{diam}(M)\ge \pi-\epsilon(n,\rho),\quad \forall\, x\in M \text{ is a $(\rho,\delta(n))$-LRRP},$$
then $\op{vol}(B_\rho(x))\ge v(n,\rho)>0$.
\endproclaim

\remark{Remark \rm 1.2} (1.2.1) Theorem 1.1 does not hold, if any one of the three
conditions is dropped; e.g., the Berger's sphere with normalizing diameter $=\pi$ ($\min \op{Ric}\to 0$), the sequence of lens spaces, $S^3/\Bbb Z_h (h\to \infty)$ ($\op{diam}(S^3_1/\Bbb Z_h)\le \frac \pi2$), or a smoothing off the two vertices of a spherical suspension of a small round sphere $S^{n-1}_\epsilon$ (given $\rho>0$, choosing $\epsilon<<\rho \delta(n)$
the two `vertices' are not $(\rho,\delta(n))$-LRRP).

\noindent (1.2.2) Theorem 1.1 may hold if one weakens the condition that all $x\in M$ are $(\rho,\delta(n))$-LRRP to that $\op{vol}(B_\rho(\tilde x))\ge v>0$. Nevertheless,
Theorem A is false under this weak local bounded covering geometry condition (\cite{An1} and \cite{Ot}).
\endremark

\proclaim{Corollary 1.3} Let $M$ be as in Theorem 1.1. Then all points in $M$ are $(c(n,\rho),\delta(n))$-LRRP, where $c(n,\rho)=\frac{\rho v(n,\rho)\op{vol}(\underline B^n_{\frac \rho2}(1))}{8\op{vol}(\underline B_\rho^n(1))^2}>0$.
\endproclaim

\demo{Proof} Let $\pi_{x_i}: (\widetilde{B_\rho(x_i)},\tilde x_i)\to (B_\rho(x_i),x_i)$ denote the Riemannian universal cover such that $x_i$ is $(\rho,\delta(n))$-LRRP. Observe
that if for any $\gamma_i\in \pi_1(B_\rho(x_i))$, $d(\gamma_i(\tilde x_i),\tilde x_i)\ge 2r_0>0$, then $B_{r_0}(x_i)$ is isometric to $B_{r_0}(\tilde x_i)$, thus $x_i$ is a $(r_0,\delta(n))$-RP.

First, $\left<\gamma_i\right>(\tilde x_i)$ is not contained in $B_{\frac \rho2}(\tilde x_i)$, where $\left<\gamma_i\right>$ denotes the subgroup generated by $\gamma_i$ such that $d(\gamma_i(\tilde x_i),\tilde x_i)\le d(\gamma_i^j(\tilde x_i),\tilde x_i)$, $\gamma_i^j\ne e$; otherwise, because $B_\rho(\tilde x_i)$ is diffeomorphic and bi-H\"older close to $\b B^n_\rho(0)$, using the Euclidean metric the energy function, $e(\cdot)=\sum_{j=1}^n\frac 12\b d^2_0(\gamma_j^j(\tilde x_i),\cdot)$, achieves the minimum at unique point $\tilde z_i\in B_{\frac \rho2}(\tilde x_i)$. Consequently, $\left<\gamma_i\right>(\tilde z_i)$ concentrates around $\tilde z_i$, a contradiction.

Let $m+1$ be the first order such that $\gamma_i^{m+1}(\tilde x_i)\notin B_{\frac \rho2}(\tilde x_i)$; $\gamma_i^j(B_{\frac \rho2}(\tilde x_i)\cap D(\tilde x_i))\subseteq B_\rho(\tilde x_i)$,$1\le j\le m$. Then $m\le \frac{\op{vol}(B_\rho(\tilde x_i))}{\op{vol}(B_{\frac \rho2}(\tilde x_i)\cap D(\tilde x_i))}\le \frac{\op{vol}(\underline B_\rho^n(1))}{\op{vol}(B_{\frac \rho2}(x_i))}$, where $D(\tilde x_i)$ denotes a fundamental domain at $\tilde x_i$. Then $\op{vol}(B_{\frac \rho2}(x_i))\ge \frac{\op{vol}(\underline B^n_{\frac \rho2}(1)}{\op{vol}(\underline B^n_\rho(1))}\cdot \op{vol}(B_\rho(x_i))\ge \frac{\op{vol}(\underline B^n_{\frac \rho2}(1))}{\op{vol}(\underline B^n_\rho(1))}\cdot v(n,\rho)$ (Theorem 1.1), thus
$$d(\gamma_i(\tilde x_i),\tilde x_i)\ge \frac\rho {2(m+1)}\ge \frac \rho{4m}\ge \frac{\rho\cdot \op{vol}(B_{\frac \rho2}(x_i))}{4\op{vol}(\b B^n_\rho(1))}\ge \frac{\rho\cdot  v(n,\rho)\op{vol}(\b B_{\frac \rho2}^n(1))}{4\op{vol}(\b B_\rho^n(1))^2}=2r_0.$$
\qed\enddemo

In our proof of Theorem A, we will use the following basic results on Ricci limit spaces.

\proclaim{Theorem 1.4} Let $M_i@>\op{GH}>>X$ be a sequence of compact $n$-manifolds such that
$$\op{Ric}_{M_i}\ge -(n-1),\quad \op{diam}(M_i)\le d.$$

\noindent {\rm (1.4.1)} {\rm (Volume convergence, \cite{Co2}, \cite{CC2})} If $\op{vol}(M_i)\ge v>0$, then $\op{vol}(M_i)\to \op{vol}(X)$, and the Hausdorff measure on $X$, $\op{vol}$,
satisfies the Bishop-Gromov volume comparison.

\noindent {\rm (1.4.2)} {\rm (Bi-H\"older convergence, \cite{CC2}, \cite{CJN})} If all $x_i\in M_i$ are $(\rho,\delta(n))$-RP, then $X$ is a manifold, and for $i$ large, $M_i$ is $\Psi(\delta |n,\rho)$ bi-H\"older diffeomorphic to $X$.
\endproclaim

Let's now present a proof of Theorem A by assuming Theorem 1.1.

\demo{Proof of Theorem A by assuming Theorem 1.1}

By a standard compactness argument, one may consider a sequence of $n$-manifolds, $M_i@>\op{GH}>>X$, satisfying that $\op{Ric}_{M_i}\ge (n-1)$, all $x_i\in M_i$ are $(\rho,\delta(n))$-LRRP, and $\op{diam}(M_i)\to \pi=\op{diam}(X)$, and show that for $i$ large, $M_i$ is bi-H\"older diffeomorphic to $S^n_1$.

By Corollary 1.3, all points in $M_i$ are $(c(n,\rho),\delta(n))$-RP. By (1.4.2), $X$ is a manifold and for $i$ large $M_i$ is bi-H\"older diffeomorphic to $X$. Let $p, q\in X$,
$d(p,q)=\op{diam}(X)=\pi$. Then for $0<r<\pi$, $B_{\pi-r}(p)\cap B_r(q)=\emptyset$. Applying the Bishop-Gromov relative volume comparison on $X$ ((1.4.1)), we derive
$$\split \op{vol}(X)&\ge \op{vol}(B_{\pi-r}(p))+\op{vol}(B_r(q))\\&=\op{vol}(X)\left(\frac{\op{vol}
(B_{\pi-r}(p))}{\op{vol}(X)}+\frac{\op{vol}(B_r(q))}{\op{vol}(X)}\right)
\\& \ge\op{vol}(X)\left(\frac{\op{vol}
(\b B^n_{\pi-r}(1))}{\op{vol}(S^n_1)}+\frac{\op{vol}(\b B^n_r(1))}{\op{vol}(S^n_1)}\right)=\op{vol}(X),\endsplit$$
thus each inequality is an equality i.e. $\frac{\op{vol}(X)}{\op{vol}(S^n_1)}=\frac{\op{vol}(B_r(q))}{\op{vol}(\underline B_r^n(1))}\to 1-\Psi(\delta(n))$ as $r\to 0$. Because $\frac{\op{vol}(M_i)}{\op{vol}(S^n_1)}\to  \frac{\op{vol}(X)}{\op{vol}(S^n_1)}$ ((1.4.1)), the desired result follows from (0.1.2) when $\delta(n)$ is suitably small.

We now verify the bi-H\"older closeness via contradiction: assuming a sequence, $\delta_i\to 0$,
and for each $\delta_i$ there is a $M_i$ of $\op{Ric}_{M_i}\ge (n-1)$, $\op{diam}(M_i) =\pi-\Psi(\delta_i|n)\to \pi$, but $d_{\op{GH}}(M_i,S^n_1)\ge \eta>0$. Because $M_i@>\op{GH}>>X$ with $\op{diam}(M_i)=\pi-\Psi(\delta_i|n)\to \pi$, repeating the above argument for a fixed $\delta$ we derive that $d_{\op{GH}}(M_i,S^n_1)<\Phi(\delta|n)$, a contradiction when $\delta=\delta_i$ sufficiently small.
\qed\enddemo

\demo{Proof of (0.4.2) using Theorem A} (Quantitative maximal diameter rigidity implies the maximal diameter rigidity)

Because $M$ is a compact manifold, there is $\delta(M)>0$ (depending on $\max \{|\op{sec}_M|\}$ and $\op{injrad}(M)$) such that for
any $0<\delta\le \delta(M)$, all points on $M$ are $(\rho,\delta)$-Reifenberg point. Let $\delta(M)=C\delta(n)$, where $\delta(n)$ is given in Theorem A. Consequently, for all $0<\delta<\frac{\delta(M)}C$, $M$ is $\Psi(\delta |n)$ bi-H\"older diffeomorphic to
$S^n_1$, thus $M$ is isometric to $S^n_1$.
\qed\enddemo

Observe that in Theorem 1.1, $(\rho,\delta(n))$-LRRP is a local regularity
condition on collapsed manifolds; which may indicate why our proof of Theorem 1.1
involves tools from several subfields: the Cheeger-Colding-Naber theory on Ricci limit spaces, local Ricci flows (\cite{Pe3}, \cite{HW1}), the collapsing structural theory of Cheeger-Fukaya-Gromov (\cite{CFG}), and the stable collapsing studied in \cite{PRT}.

For convenience of readers, we outline our approach to Theorem 1.1, with details supplied in Sections 2-4.

Arguing by contradiction, there is a sequence of compact $n$-manifolds, $M_i$, satisfying
$$\op{Ric}_{M_i}\ge (n-1),\quad M_i@>\op{GH}>>X, \quad \forall\, x_i\in M_i \text{ is $(\rho,\delta)$-LRRP}, \quad \op{vol}(M_i)\to 0,$$
and $p_i, q_i\in M_i$, $d(p_i,q_i)=\op{diam}(M_i)\to \op{diam}(X)=d(p,q)=\pi$, $p_i\to p, q_i\to q$, and we will specify the value of $\delta$ later.

Our approach is to analyze underlying structures of $M_i$ imposed by the above conditions, using which we derive a contradiction. Precisely, we find that $M_i$ contains two compact embedded infra-nilmanifolds of positive dimensions, $F(p_i), F(q_i)$, and $\pi_1(M_i)$ is isomorphic to that of a gluing of normal disk bundles of $F(p_i)$ and $F(q_i)$,
$\pi_1(M_i)\cong \pi_1(D(F(p_i))\cup_\partial D(F(q_i)))$.

For $n\ge 4$, by Van-Kampen theorem we see that $\pi_1(D(F(p_i))\cup_\partial D(F(q_i)))\cong \pi_1(F(p_i))$ which is infinite, a contradiction to (0.4.1).

For $n=3$, $F(p_i)$ and $F(q_i)$ are circles, thus $\pi_1(D(F(p_i))\cup_\partial D(F(q_i)))$ is finite unless up to a double cover, $M_i$ is homeomorphic to $S^2\times S^1$. Observe that it
suffices to prove Theorem 1.1 for the sequence of Riemannian universal cover, $\tilde M_i$, because if $\text{vol}(\tilde M_i)\ge v>0$; following the proof of Theorem A by assuming Theorem 1.1, $\tilde M_i$ is bi-H\"older diffeomorphic to $S^3_1$, thus $\op{diam}(M_i)\le \frac \pi2$, a contradiction unless $M_i=\tilde M_i$. In this case, we (directly) show to that $\tilde M_i$
admits an almost isometric $T^k$-action ($k=3-\dim(X)=1$ or $2$) without fixed point such that the diameter of each $T^k$-orbit (extrinsic or intrinsic) are uniformly converge to zero. Passing to a subsequence, we show that $\tilde M_i@>\op{GH}>>X$ is equivalent to $(\tilde M,T^k,\tilde g_i)@>\op{proj}>>(M/T^k,d_\infty)$ i.e., all $(\tilde M_i,T^k)$ are $T^k$-conjugate to $(\tilde M,T^k)$, and the orbit projection maps are $\epsilon_i$-Gromov-Hausdorff approximation (briefly, GHA), $\epsilon_i\to 0$ (Lemma 2.7). This reduction enables us to apply the stabling collapseing result in \cite{PRT} (see Lemma 2.6) to conclude that the sequence of local Riemannian universal cover of a $\rho$-tubular neighborhood of any $T^k(x)$ splits, a contradiction to (1.4.1); because $\op{Ric}_{g_i}\ge n-1$, there is a definite gap between $\op{vol}(B_\rho(\tilde x)),\tilde g_i)$ and $\op{vol}(B_\rho(\Bbb R\times \Sigma))$, where $\Sigma$ is a surface of non-negative curvature.

Finally, we explain how to get an isomorphism, $\pi_1(M_i)=\pi_1(D(F(p_i))\cup_\partial D(F(q_i)))$.

Step 1. Fixing $0<\alpha<10^{-1}$, we construct a smooth approximation to the distance function, $u_i: M_i\to \Bbb R$, $|u_i-d_{p_i}|\to 0$, such that for $i$ large $u_i$ is non-degenerate on $W_i= u_i^{-1}([\alpha,\pi-\alpha])$ (see Lemma 2.2). Consequently, $M_i$ is homeomorphic to a gluing of $U(p_i)$ and $U(q_i)$ along boundaries, $U(p_i)\cup U(q_i)=M_i-W_i$.

Observe that $\op{Ric}_{M_i}\ge n-1$ and $\op{diam}(M_i)\to \pi$ implies that the excess function, $e_{p_i,q_i}(x_i)=d(p_i,x_i)+d(q_i,x_i)-d(p_i,q_i)\to 0$ for all $x_i\in M_i$ (Lemma 2.1). This property allows one to construct, for $t_0\in (\alpha,\pi-\alpha)$, a $(1,\delta)$-splitting maps (\cite{CC1}) on a tubular neighborhood of $d_{p_i}^{-1}(t_0)$ that approximates $d_{p_i}$ (Lemma 2.2). We then obtain an open cover for $M_i-(B_{\frac \alpha2}(p_i)\cup B_{\frac \alpha2}(q_i))$, consisting of thin annulus, and $u_i$ is obtained by gluing these $(1,\delta)$-splitting maps associated to
a partition of unity, thus $u_i$ approximates $d_{p_i}$. Following the observation in \cite{Hu1}
based on the canonical Reifenberg method (\cite{CJN}), the local regularity that $x_i\in M_i$ is LRRP guarantees that $u_i$ is non-degenerate.

Step 2. We shall show that $\pi_1(U(p_i))\cong \pi_1(V(p_i))$, where $V(p_i)$ is a small perturbation of $U(p_i)$ and $V(p_i)$ is a disk bundle over a submanifold which is an infra-nilmanifold $F(p_i)$ of positive dimension. If one also assumes that $|\op{sec}_{M_i}|\le C$, using a local version of a singular nilpotent structure on a collapsed manifold with bounded sectional curvature (\cite{CFG}) it is not hard to verify the above property. In general, we will employ the Perelman's pseudo-locality applied on a local Ricci flows (\cite{HW1}) to obtain a nearby metric $g_i(t)$ on $U(p_i)$ such that $|\op{sec}_{g_i(t)}|\le C$ ($C$ is a constant independent of $i$) and $|d_{p_i}-d_{p_i,g_i(t)}|_{B_\rho(U(p_i))}<\Psi(t |n,\rho)$. Then we are able to conclude the desired result using the smoothed metrics.

\remark{Remark \rm 1.5} In view of the above our approach to Theorem A, it is worth to
note its very differences, in terms of difficulties and techniques, from proofs of topological stability mentioned in the above. For instance, in the proof of Theorem 0.5 the condition
that $\op{sec}_M\ge -\kappa$ immediately converts a small excess to the non-degeneracy of the distance function $d_p$ on $M-(B_r(p)\cup B_r(q))$ for a fixed small $r>0$, thus $M$ is homeomorphic to the gluing of two closed $r$-balls, $\bar B_r(p)\cup_\partial \bar
B_r(q)$ (\cite{GS}). The difficulty in the proof of Theorem 0.5 is to show that for $r$ sufficiently small, $d_p|_{B_r(p)-\{p\}}$ is still non-degenerate (\cite{Pe2}).
\endremark

\vskip4mm

\head 2. Proof of Theorem 1.1 (I)
\endhead

\vskip4mm

In our proof of Theorem 1.1, the condition that all $x\in M$ are LRRP is required in the proofs of
Lemma 2.2 and Lemma 2.5.

We start with the following

\proclaim{Lemma 2.1}{\rm (Small excess)} Let $ M $ be a compact Riemannian $n$-manifold of $\op{Ric}_{M} \geq (n-1)$. Given $\epsilon>0$, there exists a constant, $\Psi(\epsilon|n)$, such that if $\op{diam}(M) \geq \pi-\epsilon$, then the excess function,
$e_{p,q}(x)=d(p,x)+d(q,x)-d(p,q): M\to \Bbb R$, $e_{p,q}(x)\le \Psi(\epsilon|n)$, where $d(p,q)=\op{diam}(M)$.
\endproclaim

A technique of excess estimate was introduced in \cite{AG}, and a proof of Lemma 2.1 via
volume comparison (\cite{GP}) is similar to an alternative proof of (0.4.2) via volume comparison: let $p,q\in M$ such that $d(p,q)=\op{diam}(M)=\pi$. Then for $0<r<\pi$, $B_r(p)\cap B_{\pi-r}(q)=\emptyset$, thus $\op{vol}(M)\ge \op{vol}(B_r(p))+\op{vol}(B_{\pi-r}(q))$.
Applying Bishop-Gromov's relative volume comparison to the volume ratio of corresponding
balls, one easily sees that $\op{vol}(M)=\op{vol}(S^n_1)$, thus the maximal volume rigidity implies that $M$ is isometric to $S^n_1$.

Grove-Petersen showed that if $\op{diam}(M)>\pi-\epsilon$, then $M=B_{r+\Psi(\epsilon|n)}(p)\cup B_{\pi-r+\Psi(\epsilon|n)}(q)$ (\cite{GP}).

For any $x\in U_\alpha=U-(B_{\frac \alpha2}(p)\cup B_{\frac \alpha2}(q))$, let $A(x,r)=\bar B_{d_p(x)+r}(p)-B_{d_p(x)-r}(p)$ be the $r$-annulus centered at $p$ that contains $x$.

\proclaim{Lemma 2.2} {\rm (Small excess, LRRP and the non-degeneracy of a smooth approximation to $d_p$)}
Given $n\ge 2$, $0< \alpha <10^{-1}$, and $\rho>0$, there exist constants, $\epsilon(n,\alpha,\rho), \delta_1(n)>0$, such that for any $0<\epsilon <\epsilon(n, \alpha, \rho)$, if a region $U$ in a complete $n$-manifold $M$ contains two points, $p$ and $q$, satisfies
$$\op{Ric}_U \ge -(n-1), \quad e_{p,q}(x)<\epsilon(n,\alpha,\rho),\quad \forall\,\, x\in U \text{ is $(\rho,\delta_1(n))$-LRRP},$$
then there is a smooth function, $u: U\to \Bbb R$ such that $u$ is non-degenerate on $u^{-1}([\alpha, d(p,q)-\alpha])$, and $|u - d_{p}|< \Psi(\epsilon |n,\alpha)$.
\endproclaim

\proclaim{Corollary 2.3} Given $n\ge 2$, $0< \alpha <10^{-1}$, and $\rho>0$, there exist constants, $\epsilon(n,\alpha,\rho), \delta_1(n)>0$, such that for any $0<\epsilon <\epsilon(n,\alpha,\rho)$, if $M$ is a complete $n$-manifold satisfying
$$\op{Ric}_{M} \ge (n-1), \quad \op{diam}(M)> \pi-\epsilon,\quad \forall\,\, x\in M \text{ is $(\rho,\delta_1(n))$-LRRP},$$
then there is a smooth function $u: M\to \Bbb R$ such that $u$ is non-degenerate on $u^{-1}([\alpha, \pi -\alpha])$, and $|u-d_p| < \Psi(\epsilon |n,\alpha)$.
\endproclaim

\remark{Remark 2.4} In the proof of Lemma 2.2, $u$ is constructed by gluing of local $(1,\delta)$-splitting maps (\cite{CC1}). Because the LRRP condition, the lifting of
each local $(1,\delta)$-splitting map on the local Riemannian universal cover extends
to a local $(n,\delta)$-splitting map, which, according to (\cite{CJN}), is non-degenerate,
thus the original $(1,\delta)$-splitting map is non-degenerate when all $x\in M$ are at least $(\rho,\Psi(\delta |n))$-LRRP and $\delta$ is a uniform small constant depending on $n$. By
now one may choose $\delta_1(n)$ in Lemma 2.2 as $\Psi(\delta|n)$.
\endremark

Let $U$ be an open subset of a Riemannian $n$-manifold (e.g., $U=M$). A nilpotent structure on $U$ is an $O(n)$-invariant fiber bundle on the frame bundle over $U$, $(F(U),O(n))@>\tilde f>>(Y,O(n))$,
with fiber $N$ a nilpotent manifold and affine structure group. The $O(n)$-invariance implies that the nilpotent fibration on $F(U)$ descends to a possible singular nilpotent fibration, $U@> f>>X=Y/O(n)$, with a fiber the projection of $N$. Observe that some $f$-fiber may have dimension less than $\dim(N)$ (e.g., if $z\in F(U)$, $N_z\cap (O(n)(z))$ is a submanifold (identified with a torus subgroup $e\ne T^s<O(n)$, because which acts freely on $N_z$), then $N_z$ projects to $N_z/T^s$), or $N_z\cap (O(n)(z))$ is a finite set, thus the $f$-fiber is an infra-nilmanifold that is finitely covered by $N_z$ (e.g., the frame bundle of a flat Klein bottle $K^2$ is a $T^2$-bundle over $S^1$, and a $T^2$-fiber projects to $K^2$.)

If $U$ is also simply connected, then from the homotopy exact sequence of
$N\to F(U)\to Y$ one sees that $\pi_1(N)$ is abelian, thus the singular
nilpotent fibration on $U$ coincides with the orbits of a torus $T^k$-action.

According to \cite{Fu1} and \cite{CFG}, any collapsed manifold (or a region) with bounded sectional curvature and diameter admits a nilpotent structure
such that each $f$-fiber has positive dimension and small diameter (both intrinsically and extrinsically); see Theorem 4.1.

\proclaim{Lemma 2.5} {\rm(Local smoothing and disk-bundles)} Let $M_i@>\op{GH}>>X$ be a sequence of compact $n$-manifolds such that
$$\op{Ric}_{M_i}\ge-(n-1),\quad \op{diam}(M_i)\le d,\quad \forall\, x_i\in M_i
\text{ is $(\rho, \delta_2(n))$-LRRP}, \quad \op{vol}(M_i)\to 0.$$
Then for $x\in X$, there is $\ell_x>0$, such that for $i$ large, there is an open subset
$V_i\subset M_i$, $f_i: V_i\to B_{\ell_x}(x)$ is a singular nilpotent fibration, $f_i$ is
an $\epsilon_i$-GHA ($\epsilon_i\to 0$), $V_i$ is diffeomorphic to a disk-bundle over $F_i=f_i^{-1}(x)$.
\endproclaim

We point out that $\dim(F_i)\ge 1$, because $F_i$ is a fiber of a singular nilpotent fibration
constructed on a collapsed manifold with bounded sectional curvature (see Theorem 4.1), where
$(\rho,\delta_2(n))$-LRRP condition enables us to apply local Ricci flows technique in \cite{HW1}
to get a nearby (collapsed) metric with bounded sectional curvature.

We shall first give a proof of Theorem 1.1 for $n\ge 4$, by assuming Corollary 2.3 and Lemma 2.5.

\demo{Proof of Theorem 1.1 for $n\ge 4$}

Arguing by contradiction, assuming a sequence of $n$-manifolds, $M_i@>GH>> X$,
$$\op{Ric}_{M_i} \geq (n-1), \quad \op{diam}(M_i)\to \pi,\quad \forall\, x_i\in M_i \text{ is $(\rho, \delta(n))$-LRRP}, \quad \op{vol}(M_i) \to 0.$$
where $\delta(n)=\min\{\delta_1(n), \delta_2(n)\}$, where $\delta_1(n)$ and $\delta_2(n)$ are in Corollary 2.3 and Lemma 2.5, respectively. Let $p_i, q_i\in M_i$, $d(p_i, q_i)=\op{diam}(M_i)$, $p_i\to p$ and $q_i\to q$.

By Lemma 2.5 we may assume $\eta>0$ such that for all $i$ large, there is a
singular nilpotent fibration, $f_i: V(p_i)\to B_{2\eta}(p)$,  $B_\eta(p_i)\subset V(p_i)\subset B_{3\eta}(p_i)$, and $V(p_i)$ is a normal disk $D(p_i)$ bundle over $F(p_i)=f_i^{-1}(p)$, of extrinsic radius $\approx 2\eta$. Similarly, one gets a singular nilpotent fibration, $f_i: V(q_i)\to B_{2\eta}(q)$.

Let $\alpha=\frac {2\eta}3$. By Corollary 2.3, for $i$ large we may assume that $u_i: M_i\to \Bbb R$ is non-degenerate on $W_i=u_i^{-1}([\alpha,\pi-\alpha])$. Then $M_i=V(p_i)\cup W_i\cup V(q_i)$ (not necessarily a disjoint union). Let $U(p_i)\cup U(q_i)=M_i-W_i$. Then
$\pi_1(M_i)\cong \pi_1(U(p_i)\cup _\partial U(q_i))$.

We claim that $\pi_1(M_i)\cong \pi_1(F(p_i))$, which contains a torsion free infinite nilpotent group of finite index, a contradiction to (0.4.1).

By Corollary 2.3 and Lemma 2.5, because both $u_i$ and $d_{F(p_i)}$ approximate $d_{p_i}$,
we may assume that $\partial V(p_i)\subset u_i^{-1}([\alpha, 6\alpha])\cong \partial U(p_i)\times [\alpha,6\alpha]$ and $\partial U(p_i)\subset V(p_i)-V_{\frac 12}(p_i)\cong \partial V(p_i)\times [\frac 12,1])$, where $V_{\frac 12}(p_i)\subset V(p_i)$ is a disk-bundle over $F(p_i)$ of half size. We may assume $\op{proj}_1: \partial U(p_i)\times [\alpha,6\alpha]\to \partial U(p_i)$ and $\op{proj}_2: \partial V(p_i)\times [\frac 12,1]\to \partial V(p_i)$. Then the composition of maps, $\op{proj}_2\circ \op{proj}_1\circ \op{incl}_{\partial U(p_i)}: \partial U(p_i)\to \partial V(p_i)$, induces an isomorphism on fundamental groups, thus $\pi_1(\partial U(p_i))\cong \pi_1(\partial V(p_i))$. Similarly, one gets that $\pi_1(U(p_i))\cong \pi_1(V(p_i))$, and thus
$\pi_1(M_i)\cong \pi_1(V(p_i)\cup_\partial V(q_i))$.

Next, we claim that if $\dim(D(p_i))=2$, then $\dim(D(q_i))=2$ and $n=3$. First, because the universal cover of $\partial (V(p_i))$ is an Euclidean space, and because $\partial V(p_i)$ is homeomorphic to $\partial V(q_i)$, $\dim(D(q_i))=2$ (if $\dim(D(q_i))\ge 3$, then universal cover of $\partial V(q_i)$ is not an Euclidean space). Secondly, if $\dim(F(p_i))=\dim(F(q_i))\ge 2$, by Van Kampen theorem it is easy to see that $\pi_1(M_i)\cong \pi_1(V(p_i)\cup_\partial  V(q_i))$ is infinite (note that $D(q_i)$ (resp. $D(p_i)$) kills only one homotopy class in $F(p_i)$ (resp. $F(q_i)$), a contradiction to (0.4.1)). By now the claim has been verified.

In view of the above, $n\ge 4$ implies that $\dim(D(p_i))\ge 3$ i.e., $\partial D(p_i)$ is a sphere of dimension $\ge 2$. Because $\pi_1(M_i)\cong \pi_1(V(p_i)\cup_\partial  V(q_i))$, by Van Kampen theorem $\pi_1(M_i)\cong [\pi_1(V(p_i))*\pi_1(V(q_i))]/\pi_1(\partial V(p_i))\sim \pi_1(\partial V(q_i))$. Because $\pi_1(\partial V(p_i))\cong \pi_1(V(p_i))\cong \pi_1(F(p_i))$ (the same holds for $q_i$) and $\pi_1(\partial V(p_i))\cong \pi_1(\partial V(q_i))$, $\pi_1(M_i)\cong \pi_1(F(p_i))$, a contradiction to (0.4.1).
\qed\enddemo

Because for $n=3$, $F(p_i)$ (resp. $F(q_i)$) is a circle and $D(p_i)$ is a $2$-disk, $\pi_1(M_i)$ is finite excepts a double cover of $M_i$ is diffeomorphic to $S^2\times S^1$;
so $\pi_1(M_i)$ is generally a finite group (no contradiction to (0.4.1)).
Instead, we shall employ the splitting property of a stable collapsing (\cite{PRT})
to derive a contradiction to (1.4.1).

The notion of a stable collapsing was introduced in \cite{PRT}. In the following,
we will focus on a special stable collapsing as follows: a collapsing sequence of compact $n$-manifolds, $M_i@>\op{GH}>>X$, is called a stable collapsing,
if $M_i=(M,g_i)$, $M$ admits a $T^k$-action without fixed points, such that $X$ is homeomorphic to $M/T^k$ and orbit projection, $\op{proj}: (M,g_i)\to M/T^k$ (equipped with the metric on $X$), is an $\epsilon_i$-Gromov-Hausdorff approximation (briefly, GHA) with $\epsilon_i\to 0$; thus the intrinsic diameter of all $T^k$-orbits uniformly converge to zero while the metric orthogonal to
$T^k$-orbits converges.

A simple example of a stable collapsing is the Berger's sphere,
an one-parameter family of $S^1$-invariant metrics on $S^3_1$, collapsing to its orbit space of the Hopf-fibration $S^3_1$, which is a round $2$-sphere of radius $\frac 12$.

Geometric properties of a stable collapsing was studied in \cite{PRT}; following the proof of Theorem 0.3 in \cite{PRT}, one gets the following result:

\proclaim{Lemma 2.6} {\rm (Stabling collapsing, local splitting)} Let
$(M,T^k, g_i)$ be a stable collapsing which satisfies
$$\op{Ric}_{g_i}\ge -(n-1)\delta_i\to 0, \quad \op{diam}(g_i)\le d, \quad \forall\, x\in M \text{ is $(1,\delta(n))$-LRRP w.r.t. all $g_i$}.$$
Then for any $x\in M$, $(\widetilde {B_{\delta(n)}(x)},\tilde x,\tilde g_i)@>\op{GH}>> (Z\times \Bbb R^\ell,x)$, for some $1\le \ell\le k$.
\endproclaim

\demo{Proof} Because Lemma 2.6 is easily seen follows from the proof of Theorem 0.3 in \cite{PRT}, we will only point out two key observations in the proof.

Note that a stable collapsing sequence of metrics $g_i$ on $M$ in Theorem 0.3 satisfy that $\lambda\le \op{sec}(g_i)\le \Lambda$, which implies a fixed nilpotent structure on $M$ (see Theorem 4.1) and every orbit collapses to a point while metrics orthogonal to orbits converges.

The first key observation when given a stable collapsing, $(M,T^k,g_i)@>\op{GH}>>(M/T^k,d_\infty)$,
in Lemma 2.6, the gluing operation can be performed under weak curvature conditions:
e.g., $\op{Ric}_{M_i}\ge -(n-1)$ and $(\rho,\delta(n))$-LRRP (which holds under $\lambda\le \op{sec}(g_i)\le \Lambda$). Precisely, given an atlas, $\{(U_\alpha,T^k)\}$, $U_\alpha$ is
a tubular neighborhood of $T^k(x_\alpha)$ (the slice theorem), the gluing operation
requires an equivariant GH-convergence, $(\tilde U_\alpha,A^k,\tilde g_i)$, and the topological stability through the GH-convergence, where $A^k$ denotes the lifting group of $T^k$ that
(effectively) acts on the universal cover $\tilde U_\alpha$ of $U_\alpha$ (e.g., $A^k=T^s\times \Bbb R^{k-s}$ with $T^s$ the isotropy group at $x_\beta$). Note that the complete GH-limit is from a sequence of $\tilde g_i$ on a region in $\tilde U_\alpha$, consisting of points away from $\partial \tilde U_\alpha$ (\cite{Xu}).

The other key observation is that the condition, $\op{Ric}_{M_i}\ge -\delta_i(n-1)\to 0$, implies that fixing a chart $U_{\beta}$, when a complete GH-limit of a sequence $\tilde g_i$ on $\tilde U_\beta$ contains a line, the complete GH-limit space splits off a line (\cite{CC1}). Consequently, for every $U_\alpha$, the corresponding complete GH-limit splits off a line (e.g., if $U_\beta\cap U_\alpha\ne \emptyset$, then the intersection contains a line).

A desired chart (or a desired atlas) can be constructed as follows: let
$x_\beta\in M$ such that either $T^k(x_\beta)$ has a maximal isotropy group (if the $T^k$-action is not free), or (if the $T^k$-action is free) the maximum of the area density on $T^k(x_\beta)$ (with induced metric) is minimum among all $x\in M$ (e.g., if $T^k$ acts
by isometries, then the area density function is a constant on each $T^k$-orbit i.e. a function defined on $M/T^k$). It is easy to see that $(\pi_\beta^{-1}(T^k(x_\beta)),\tilde g_i)$
converges to an isometrically embedded Euclidean space in the the complete GH-limit (e.g., if $T^k$ acts isometrically, then $\pi^{-1}_\beta(T^k(x_\beta))$ is a totally geodesic Euclidean space, so is its GH-limit `totally geodesic' and Euclidean space. In general, the limit group of $A^k$ (which is isomorphic to $A^k$) acts isometrically on the complete GH-limit.
\qed\enddemo

In our proof of Theorem 1.1 for $n=3$, the main technical result is the following.

\proclaim{Lemma 2.7} Let $M_i@>\op{GH}>>X$ be a sequence of compact $3$-manifolds satisfying
$$\op{Ric}_{M_i}\ge 2,\quad \op{diam}(M_i)\to \pi,\quad \forall\, x_i\in M_i \text{ is $(\rho,\delta(3))$-LRRP},\quad \op{vol}(M_i)\to 0.$$
If $\pi_1(M_i)=0$, then $M_i$ contains a subsequence, $M_i=(M,g_i)$, such that $M$ admits
a $T^k$-action without fixed point ($k=\dim(X)$), and the orbit projection map, $\op{proj}: (M,T^k,g_i)\to (M/T^k,d_\infty)$, is an $\epsilon_i$-GHA, $\epsilon_i\to 0$.
\endproclaim

\demo{Proof} We first construct $T^k$-action on $(M_i,g_i)$ without fixed point, which almost preserves $g_i$.

As seen in the proof of Theorem 1.1 (I) and (II), $M_i=U(p_i)\cup W_i\cup U(q_i)$. For $n=3$,
$X=[0,\pi]$ or $X=S(S^1_r)$, a spherical suspension over a circle of radius $r\le 1$ (note that $V(p_i)$ is a disk bundle implies that if $\dim(X)=2$, then $\partial X=\emptyset$ and $X=S(S^1_r)$). Passing to a subsequence and for points always from the boundaries, we may assume local Riemannian universal covers,
$$\CD (\widetilde{U(p_i)},\tilde p_i,\Gamma_i)@>\op{eqGH}>>(\tilde X,\tilde p,G)\\
@V \pi_i VV @V\op{proj} VV\\
(U(p_i),p_i)@>\op{GH}>>(B_\rho(p),p),
\endCD$$
where $\Gamma_i=\pi_1(U(p_i))$.

Case 1. $X=[0,\pi]$. Then $W_i$ is homeomorphic to $T^2\times [\alpha,\pi-\alpha]$, with
each $T^2$-fiber uniformly converges to a point, and $U(p_i)$ is homeomorphic to
$S^1\times D^2$ ($\pi_1(M_i)=0$). Clearly, the $T^2$-action on $\partial U(p_i)$
extends to a $T^2$-action on $U(p_i)$ and on $U(q_i)$, without fixed point, and
the $T^2$ acts almost isometrically such that every $T^2$-orbit collapses to a point.
Because $M_i$ is simply connected and $M_i$ is a gluing of two solid torus along
boundaries, we identify $M_i$ to $(S^3,T^2)$ and $M_i\to X$ is equivalent to a stable collapsing, $(S^3,T^2,g_i)\to S^3/T^2,d_\infty)$.

Case 2. $X=S(S^1_r)$. Then $W_i$ is homeomorphic to $S^1_{\epsilon_i}\times (S^1\times [\alpha,\pi-\alpha])$. It is clear that the $S^1_{\epsilon_i}$-action uniquely extends to an $S^1$ on $U(p_i)$ and $U(q_i)$. Because the $S^1$-orbit has a small diameter and the LRRP condition, the $S^1$-action has no fixed point. Again, we now identify $M_i$ as $S^3$, and the conjugate class of the $S^1$-action (constructed from $g_i$) is determined by the isotropy groups at $p_i$ and $q_i$, whose order
is uniformly bounded above, passing to a subsequence we may assume that all $S^1$-actions on $S^3$ are conjugate. By now it is clear that $M_i@>\op{GH}>>X$ is equivalent to a stable collapsing,
$(S^3,S^1,g_i)@>\op{GH}>>(S^3/S^1,d_\infty)$.
\qed\enddemo

\demo{Proof of Theorem 1.1 for $n=3$}

Arguing by contradiction, assuming a sequence of compact $3$-manifolds, $M_i@>\op{GH}>>X$, satisfying
$$\op{Ric}_{M_i}\ge 2,\,d(p_i,q_i)=\op{diam}(M_i)\to \pi,\,\forall\, x_i\in M_i \text{ is $(\rho,\delta(3))$-LRRP},\, \op{vol}(M_i)\to 0,$$
and $p_i\to p,, q_i\to  q$.

We first assume that $M_i$ is simply connected. By Lemma 2.7, passing to a subsequence we obtain a stable collapsing, $(S^3,T^k,g_i)@>\op{GH}>>(S^3/T^k,d_\infty)$, such that $\op{proj}: (S^3,T^k,g_i)\to (S^3/T^k,d_\infty)$ is an $\epsilon_i$-GHA, $\epsilon_i\to 0$.

Because $X=S^3/T^k$ is a spherical suspension of dimension $1$ or $2$, there is $x\in X$ and
$\rho_0=\rho(X)>0$, such that $B_{\rho_0}(x)$ is isometric to an interval of length $2r_0$, or to a $\rho_0$-ball in $S^2_1$. Let $x_i\in (M,g_i)$, $x_i\to x$. By Lemma 2.6, $(\pi^{-1}(B_{\rho_0}(x_i)),\tilde g_i)@>\op{GH}>> (x-\rho_0,x+\rho_0)\times \Bbb R^2$, or $\b B^2_{\rho_0}(1)\times \Bbb R^1$. By volume comparison and (1.4.1), we derive a contradiction:
$$\split \op{vol}(B_{\rho_0}(\tilde x))&=\lim_{i\to \infty}\op{vol}(B_{\rho_0}(\tilde x_i,\tilde g_i))\le \op{vol}(\b B^3_{\rho_0}(1))\\& <
\min\{\op{vol}(B_{\rho_0}(\frac \pi2, [0,\pi]\times \Bbb R^2)),\op{vol}(B_{\rho_0}(S^2_1\times \Bbb R^1))\} \\& \le \op{vol}(B_{\rho_0}(\tilde x)).\endsplit$$
\qed\enddemo

\vskip4mm

\head 3. Proof of Theorem 1.1 (II)
\endhead

\vskip4mm

In this section, our goal is to prove Lemma 2.2.

\proclaim{Lemma 3.1} Let the assumptions be as in Lemma 2.2. For any $x\in U_\alpha$, there is $r_x\ge \Psi(\alpha|n)>0$, such that if $e_{p,q}(x)|_{U_\alpha}<\epsilon$, then the
solution
$$\cases \Delta u=0 & A(x,2r_x)\\ u=d_p & \partial A(x,2r_x)\endcases$$
satisfies that $u$ is non-degenerate on $A(x,r_x)$ and $|u-d_p|_{A(x,r_x)}<\Psi(\epsilon | n,\alpha)$.
\endproclaim

\demo{Proof} For $x\in U_\alpha$, $L=\min \{d_p(x), d_q(x)\}\ge \frac \alpha2$. Because $e_{p,q}|_{U_\alpha}<\epsilon<<L$, there is $r_x\ge \Psi(\alpha|n)>0$ (note that $\alpha$ is a fixed constant) such that $u$ is a $(1,\delta)$-splitting
map satisfying the following estimates (\cite{CC1}, \cite{Ch}):

\noindent (3.1.1) $|u-d_p|_{A(x,2r_x)}<\Psi(\epsilon |n,\alpha)$,

\noindent (3.1.2) $|\nabla u|_{A(x,2r_x)}\le 1+\Psi(\epsilon|n,\alpha)$;

\noindent (3.1.3) For any $y\in A(x,r_x)$, $\Psi(\alpha|n)\le r\le r_x$,
$$\frac 1{\op{vol}(B_r(y))}\int_{B_r(y)}||\nabla u|^2-1|d\op{vol}<\Psi(\epsilon|n,\alpha).$$


The non-degeneracy of $u$ follows from the observation in \cite{Hu1}: because $x\in A(x, r_x)\cap U_\alpha$ is $(\rho,\delta_1(n))$-LRRP, the lifting of $u$ over $\pi^{-1}(B_{2r_x}(x))\subset \widetilde{B_\rho(x))}$ extends to a $(n,\delta)$-splitting map which is non-degenerate (\cite{CJN}), provide that $\delta=\delta_1(n)$ small. In particular, the lifting of $u$ is non-degenerate
i.e., $u$ is non-degenerate.
\qed\enddemo

Because $U_\alpha$ is compact, using Lemma 3.1 we may construct an open cover for $U_\alpha$
by a finite number of thin annulus, $\{(A(x_i,r_i),u_i)\}$,  such that $u_i$ is a non-degenerated harmonic approximation of $d_p|_{A(x_i,r_i)}$. Moreover, $\{x_i\}$ can be chosen so that $d_p(x_j)>d_p(x_{j+1})$, and $A(x_j,r_j)\cap A(x_k,r_k)\ne \emptyset $ if and only if $k=i-1$
or $j+1$.

Let $\{f_i\}$ be a partition of unity associate to $\{A(x_i,r_i)\}$ such that $|\nabla f_i|, |\Delta f_i|\le C(n,\alpha)$ (\cite{CC1}), and let
$$u=\sum_if_iu_i,\quad x\in U_\alpha.$$

\proclaim{Lemma 3.2} Let the assumptions be as in Lemma 2.2, and let $u$ be as in the
above. Then $u$ is non-degenerate on $U_\alpha$ and $|u-d_p|<\Psi(\epsilon,\delta|n,\alpha)$.
\endproclaim

Note that Lemma 2.2 follows from Lemma 2.2 and Lemma 3.2. In the proof of Lemma 3.2, we shall use the following reduced criterion on the non-degenerace.

\proclaim{Lemma 3.3}{\rm (A reduced criterion of non-degeneracy \cite{Hu2})} Given $n\ge2, 1\ge \rho>0$, and $C>0$, there exists a constant, $\delta(n,C)>0$, $\rho\ge 2$, such that if a complete $n$-manifold $M$ satisfies
$$\op{Ric}_{M} \ge -(n-1)\delta(n,C), \quad x \text{ is a $(\rho, \delta(n,C))$-LRRP},$$
and if a map $ u: B_{2}(x)\to \Bbb R^{k}$ satisfies the following properties:


\noindent{\rm (3.3.1)} $|\Delta u |\le C $,

\noindent{\rm (3.3.2)} $|\nabla u|\le 1 +\delta(n,C)$,
	
\noindent{\rm (3.3.3)} $\frac{1}{\op{vol}(B_{2}(x))}\int_{B_{2}(x)} |
|\nabla u|^{2}-1||d\op{vol}\le \delta(n,C)$,	

\noindent then $u$ is non-degenerate on $B_{1}(x)$.
\endproclaim

The non-degeneracy of $u$ is a local property, thus the non-degeneracy holds if one replaces ``$B_2(x)$'' by ``$B_r(x)$'', $0<r<1$. Note that a function $u$ satisfying (3.3.1)-(3.3.3) is also called a (generalized) $(1,\delta)$-splitting map.

\demo{Proof of Lemma 3.2}

Let $u=\sum_jf_ju_j$ be defined in the above; $u_j$ (3.1.1)-(3.1.4) with $\Psi_1(\epsilon|n,\alpha)$, and
$|\nabla f_j|, |\Delta f_j|\le C(n,\alpha)$.

First, by  (3.1.1) we derive
$$\split
		|u(x) - d_p(x)| & = |\sum_{j=1}^{N}f_j(x)u_j(x) - \sum_{j=1}^{N}f_jd_p(x)| \\
		& = \sum_{j=1}^{N} f_j |u_j(x) - d_p(x)|\le \sum_{j=1}^N f_j\Psi_1(\epsilon|n,\alpha) \\
		& =\Psi_1(\epsilon| n,\alpha).\endsplit \tag 3.3.4$$

Following a standard argument presented below (cf. \cite{Ch}), (3.3.4) is all that is required to verify (3.3.1)-(3.3.3) i.e., $u$ is non-degenerate (Lemma 3.3).

Since that each $u_j$ satisfies (3.1.1) and $x\in A(x_j,r_j)\cap A(x_k,r_k)\ne \emptyset$ if and only
if $k=j+1$ or $j-1$,
$$\split |u_j(x) - u_k(x)| &\le |u_j(x) -d_{p}(x)| + |u_k(x) - d_{p}(x)| \\&\le 2\Psi_1(\epsilon|n, \alpha)=\Psi_2(\epsilon|n,\alpha).\endsplit $$
Note that $u_j- u_k$ is a harmonic function on $A(x_j,2r_j)\cap A(x_k,2r_k)$, by Cheng-Yau gradient estimate (\cite{CY})
$$\sup_{A(x,r_j)\cap A(x,r_{j+1})}\{|\nabla u_j-\nabla u_k|\}\le C(n, \alpha)\Psi(\epsilon|n,\alpha)=\Psi_3(\epsilon|n,\alpha). $$

For any $x\in U_\alpha$, we now verify (3.3.1) as follows:
$$\split
		| \Delta u(x) | & =  |\Delta (f_ju_j+f_ku_k)-\Delta u_k|=|\Delta(f_ju_j-f_ju_k)|\\&=|\Delta f_j|\cdot |u_j-u_k|+2|\nabla f_j|\cdot |\nabla (u_j-u_k)|\\&\le C(n,\alpha)\Psi_2(\epsilon|n,\alpha)+2C(n,\alpha)\Psi_3(\epsilon|n,\alpha)=\Psi_4(\epsilon|n,\alpha).
\endsplit$$
Note that we shall verify (3.3.1)-(3.3.3) in Lemma 3.3 with $C=\Psi_4(\epsilon|n)$, without loss of generality we may assume
$\delta(n,C)\le 1$.
Because
$$\split |\nabla u - \nabla u_j|&=|\nabla f_j(u_j-u_k)|\\&=|\nabla f_j|\cdot |u_j-u_k|+|f_j|\cdot |\nabla (u_j-u_k)|\\& \le C(n, \alpha)\Psi_2(\epsilon|n, \alpha)+\Psi_2(\epsilon|n,\alpha)=\Psi_5(\epsilon|n,\alpha). \endsplit \tag 3.3.5$$
by (3.1.2) we verify (3.3.2):
$$|\nabla u|\le |\nabla u - \nabla u_{j}| + |\nabla u_j|\le \Psi_5(\epsilon|n,\alpha)+1+\Psi_1(\epsilon|n)<1+\Psi_6(\epsilon|n,\alpha).$$
Consequently,
$$|\nabla u+\nabla u_j|\le |\nabla u|+|\nabla u_j|\le 3.\tag 3.3.6$$
Finally, we verify (3.3.3) below:
$$\split \frac{1}{\op{vol}(B_{\Psi(\alpha|n)}(x))}&\int_{B_{\Psi(\alpha|n)}(x)}  |  | \nabla u |^{2} -1 |d\op{vol}
 \\ &\le \frac{1}{\op{vol}(B_{\Psi(\alpha|n)}(x))}\int_{B_{\Psi(\alpha|n)}(x)}(||\nabla u_i |^2-1|+||\nabla u|^2-|\nabla u_i|^2|)d\op{vol}\\&\le \frac{1}{\op{vol}(B_{\Psi(\alpha|n)}(x))}\int_{B_{\Psi(\alpha|n)}(x)}(| |\nabla u_i|^2-1|+ |\nabla u-\nabla u_i||\nabla u+\nabla u_i|)d\op{vol} \\& \le \Psi_1(\epsilon|n,\alpha)+3\Psi_5(\epsilon|n,\alpha). \hskip6mm \text{(by (3.3.3), (3.3.5) and (3.3.6))}\endsplit$$
\qed\enddemo

\vskip4mm

\head 4. Proof of Theorem 1.1 (III)
\endhead

\vskip4mm

In this section, we will prove Lemma 2.5, thus we complete the proof of Theorem 1.1.

We first prove Lemma 2.5 under a strong condition: $|\op{sec}_{M_i}|\le C$, a constant for all
$i$. In this case, we will apply a local version of the following singular nilpotent fibration theorem.

\proclaim{Theorem 4.1}{\rm (Singular nilpotent fibration, \cite{Fu1-2}, \cite{CFG})}
Let $M_i@>\op{GH}>>X$ be a sequence of complete $n$-manifolds satisfying
$$|\op{sec}_{M_i}|\le 1, \quad \op{diam}(M_i) \leq d, \quad \op{vol}(M_i) \to 0.$$
Then the sequence of the frame bundles equipped with a canonical metrics contains
a subsequence equivariant converges, $(F(M_i),O(n))@>\op{eqGH}>>(Y,O(n))$, such that
$Y$ is a $C^{1,\alpha}$-manifold, and for $i$ large, there is
an $O(n)$-invariant fiber bundle map, $\tilde f_i: (F(M_i),O(n))\to (Y,O(n))$, satisfying

\noindent{\rm (4.1.1)} $\tilde f_i$ is $\epsilon_i$-GHA, $\epsilon_i=d_{\op{GH}}(F(M_i),Y)\to 0$.
		
\noindent{\rm (4.1.2)} $\tilde f_i$ is a $\Psi(\epsilon_i)$ Riemannian submersion i.e., for any
vector orthogonal to $\tilde f_i$-fiber, $e^{-\Psi(\epsilon_i)} \le \frac{|df_i(v_{i})|}{|v_i|} \le e^{\Psi(\epsilon_i)}$.

\noindent{\rm (4.1.3)} The fundamental form of each fiber $\tilde f_i$-fiber,
$|\op{II}_{\tilde f_i}| \le C(n).$
		
\noindent{\rm (4.1.4)} A $\tilde f_i$-fiber is diffeomorphic to a nilmanifold whose projection
to $M_i$ has positive dimension, and the structural group reduces to affine.
\endproclaim

By the $O(n)$-invariance, $\tilde f_i$ descends to a singular nilpotent fibration and  $\epsilon_i$-GHA, $f: M_i\to X=Y/O(n)$, each $f_i$ is an infra-nilmanifold, and the following diagram commutes:
$$\CD (F(M_{i}), O(n)) @>\tilde f_i>> (Y, O(n)),\\
@V\op{proj}_i VV  @V \op{proj} VV\\
M_i@> f_i>>X=Y/O(n),
\endCD\tag 4.1.5$$

Because $\op{diam}(f_i^{-1}(x))$ are uniformly small, $f^{-1}_i(B_r(x))$ approximates to $B_r(x_i)$ ($f_i(x_i)=x$), and because every $f_i$-fiber is an infra-nilmanifold of
positive dimension, $M_i$ satisfies a local bounded covering geometry. On the other hand,
a non-collapsed $M$ trivially satisfies a local bounded covering geometry, we conclude
the following property.

\proclaim{Lemma 4.2} {\rm (Bounded covering geometry, \cite{CFG})} Given $n\ge 2$, there exist constants, $\rho(n), \epsilon(n)>0$, such that
if $M$ is a complete $n$-manifold of $|\op{sec}_M|\le 1$, then for any $x\in M$, there is
a neighborhood $U\supset B_{\rho(n)}(x)$, $\pi: (\tilde U,\tilde x)\to (U,x)$ is the Riemannian universal cover, then
$\op{injrad}(B_{\rho(n)}(\tilde x))\ge \epsilon(n)$.
\endproclaim

\demo{Proof of Lemma 2.5 for the case that $|\op{sec}_{M_i}|\le C$}

For any $x\in X$, we shall find $\ell_x>0$ such that $V_i=f_i^{-1}(B_{\ell_x}(x))$ is diffeomorphic to
a disk bundle over $F_i$, an infra-nilmanifold.

By Theorem 4.1, for $i$ large we may assume a singular nilpotent fibration, $f_i: M_i\to X$.
For $x\in X$, passing to a subsequence, by Lemma 4.2 we may assume the following commutative equivariant convergence with regularity:
$$\CD (\widetilde {f_i^{-1}(B_{\frac {\rho(n)}2}(x))},\tilde x_i,\Gamma_i)@>C^{1,\alpha}-\op{eqGH}>>(\tilde X,\tilde x,G)\\
@V\pi_i VV @V\op{proj} VV\\
(f_i^{-1}(B_{\frac {\rho(n)}2}(x)),x_i)@>f_i>>(B_{\frac{\rho(n)}2}(x)=\tilde X/G,x),
\endCD$$
where $f_i(x_i)=x$, $\pi_i: (\widetilde{f_i^{-1}(B_{\rho(n)}(x))},\tilde x_i)\to (f_i^{-1}(B_{\rho(n)}(x)),x_i)$ is the Riemannian universal cover, $\Gamma_i=\op{im}[\pi_1(f_i^{-1}(B_{\frac {\rho(n)}2}(x)),x_i)\to \pi_1(B_{\rho(n)}(x_i),x_i)]$, and the limit group $G$ is a closed Lie group. Without loss of generality, we may assume that $B_{\rho(n)}(x)$ is contractible. Then the singular nilpotent structure on $f_i^{-1}(B_{\frac {\rho(n)}2}(x))$ has the following alternative formulation: there is a connected nilpotent Lie group $G_0$ that acts isometrically on
$\widetilde{f^{-1}_i(B_{\rho(n)}(x))}$, that extends to an isometric $\Gamma_i\rtimes G_0$-action. The following properties hold:

\noindent (2.5.1) $\pi_1(f_i^{-1}(B_{\frac {\rho(n)}2}(x)),x_i)\cong \pi_1(F_i,x_i)$, thus $\Gamma_i$ is torsion free and contains a nilpotent subgroup of bounded index.

\noindent (2.5.2) The $\Gamma_i$-action on $\widetilde{f_i^{-1}(B_{\rho(n)}(x))}$ preserves every $G_0$-orbit, $\pi_i(G_0(\tilde z_i))$ is a $f_i$-fiber.

\noindent (2.5.3) For $i$ large, $G_0$ is isomorphic to the identity component of $G$ acting
on $\tilde X$.

Note that $G_0$ is a nilpotent Lie group acting isometrically on $\tilde X$. Then the $G_0$-action has a normal slice at $\tilde x$ ($\op{proj}(\tilde x)=x$) (\cite{Pa}), $D_{2\ell_x}(\tilde x)$, of radius $2\ell_x>0$. By the (pointed) $C^{1,\alpha}$ equivariant convergence it is clear that
for $i$ large the $G_0$-action on $\pi_i^{-1}(B_{\ell_x}(F_i))$ has a normal slice at $\tilde x_i$ of radius $\ell_x$. By (2.5.3), for any $\gamma_i\in \Gamma_i$, $\gamma_i(D_{\ell_x}(\tilde x_i))\cap D_{\ell_x}(\tilde x_i)=\emptyset$, thus the normal slices along $\pi_i^{-1}(F_i)$ descends to the desired normal disk bundle along $F_i$.
\qed\enddemo

In view of the above proof, we see that key ingredients in the above proof are a singular nilpotent fibration structure, a corresponding slice lemma and
regularity of equivariant convergence on local universal cover. Our approach to Lemma 2.5 in
general is to show that locally the metric can be approximated by one with bounded sectional curvature, so we able to apply the proof for the case of bounded sectional curvature.
Here the local smoothing tool is from \cite{HW1} on local Ricci flows.

The Ricci flows $g(t)$ on a Riemannian manifold $(M,g)$ was introduced by Hamilton (\cite{Ha}),
where $g(t)$ are the solutions of the following parabolic equation:
$$\frac {\partial} {\partial t}(g(t))=-2\op{Ric}(g(t)),\quad g(0)=g.$$
Since then, Ricci flows (among others things) has been a powerful tool in improving regularities of a metric on $M$. Because a Ricci flows distincts when $|\op{Rm}(g(t))|$ blows up, an estimate on $|\nabla \op{Rm}(g(t))|$ implies a definite flow time, $t_0>0$, that, in turn, gives a uniform bound on $|\op{Rm}(g(t_0))|$. In \cite{Sh}, an estimate on $|\nabla^k\op{Rm}(g(t))|$ was obtained
under a strong curvature condition on $g$, $|\op{sec}_g|\le C$ (here $\op{vol}(B_1(x))$ may be
very small). In \cite{Pe3}, an estimate on $|\op{Rm}(g(t))|$ was obtained (referred as the Perelman's pseudo-locality) under a weak curvature condition on $g$, the scalar curvature $S(g)\ge -C$, and a local Euclidean isoperimetric inequality condition (here $\op{vol}(B_1(x))$ is almost equal to $\op{vol}(B_1^n(0))$). Moreover, a manifold of $\op{Ric}_M\ge -(n-1)$ and every $x\in M$ is $(\rho,\delta(n))$-RP satisfies the Euclidean isoperimetric inequality condition (\cite{CM}, \cite{CRX1}).

In \cite{HW1}, a local version of the Perel'man's pseudo-locality method was obtained on a manifold $M$ of $\op{Ric}_M\ge -(n-1)$ and any $x$ is a $(\rho,\delta(n))$-LRRP (here $\op{vol}(B_1(x))$ may be very small but $\op{vol}(B_1(\tilde x))$ almost equals to
$\op{vol}(\b B_1^n)$).

Roughly, the $(\rho,\delta(n))$-LRRP condition allows one to conformally extend $g|_{B_{\frac \rho2}(x)}$ to a complete metric $\hat g$ on $B_\rho(x)$, such that $|\op{sec}_{\hat g}|\le C<\infty$ and all points in the Riemannian universal covering of $(B_\rho(x),\hat g)$ are $(\rho,\delta(n))$-RP. Then one applies the Perel'man's pseudo-locality method to conclude a desired estimate on $|\op{Rm}(g(t))|$ (the Ricci flows on $\widetilde{B_\rho(x)}$ is the pullback of the Ricci flows on $B_\rho(x)$, starting with $\hat g$).

In view of the above, we may state a result in \cite{HW1} in the following form for our purpose.

\proclaim{Theorem 4.3}{\rm (Local Ricci flows \cite{HW1})}
Given $n \ge 2$, $\rho>0$ and $\alpha_{0}\in (0, 10^{-1})$, there exist constants, $1>\delta(n, \alpha_{0}), \epsilon(n,\rho,\alpha_{0})>0$,  such that if a complete $n$-manifold $(M, g)$ satisfies $\op{Ric}_{g}\ge -(n-1)$ and a point $x\in M$ is $(\rho, \delta_{2}(n, \alpha_{0}))$-LRRP, then there exists a Ricci flow $g(t)$ on $B_{\frac \rho2}(x,g)$ with $t\in (0, \epsilon^{2}(n,\rho,\alpha_{0})]$ such that
	$$\sup_{B_{\frac \rho2}(x, g)}| \op{Rm}(g(t)) | \leq \frac{\alpha_{0}}{t} + \frac{1}{\epsilon^{2}(n,\rho,\alpha_{0})}.$$
\endproclaim

To make sure the smoothing metric by local Ricci flows contains essentially same collapsing geometry of
$g(0)$, one has to show that the smoothing metric
is close to the original metric i.e. the following distance estimate.

\proclaim{Lemma 4.4}{\rm (Distance distortion estimate \cite{CRX1}, \cite{BW}, \cite{HW1}, \cite{HKRX})}
Let the assumption be as Theorem 4.3.

\noindent {\rm (4.4.1)} For any $x, y \in B_{\frac \rho2}(x, g)$ with $d_{g(t)}(x, y) \leq \sqrt{t}$,
$$|d_{g(0)}(x, y) - d_{g(t)}(x, y)|\le \Psi(\alpha_{0}|n)\sqrt{t}.$$

\noindent {\rm (4.4.2)} For any $x, y \in B_{\frac \rho2}(x, g)$ with $d_{g(0)}(x, y)\le 1$,
$$c(n)^{-1}d_{g(t)}(x, y)^{1+\Psi(\alpha_{0}|n)}\le d_{g(0)}(x, y)\le c(n)d_{g(t)}(x, y)^{1-\Psi(\alpha_{0}|n)}.$$
\endproclaim


\demo{Proof of Lemma 2.5}

Let $(B_\rho(x_i), g_i) @>\op{GH}>> (B_\rho(x), d_\infty)$ be as in Lemma 2.5. By
Theorem 4.3 (with fixed small $\alpha_0$ so that $\Psi(\alpha_0|n)<10^{-4}$),
for each $g_i$ there is a local
Ricci flows $g_{i}(t)$ with $t\in (0, \epsilon(n,\rho,\alpha_{0})^{2}]$, such that

\noindent{\rm (2.5.4)} $|\sec_{g_{i}(t)}| \leq \frac {\alpha_{0}}t + \frac 1{\epsilon_{0}(n, \rho,\alpha_0)^{2}}$,

\noindent{\rm (2.5.5)} $\op{vol}_{g_{i}(t)}(B_{\frac \rho8}(x_{i})) \to 0$.

By Lemma 4.4, taking subsequence if necessary we have
$$\CD
(B_{\frac \rho4}(x_{i}), g_{i}(t)) @>\op{GH}>> (B_{\frac \rho4}(x), d_{t}) \\
@V\op{id}_{B_{\frac \rho4}(x_i)}VV @Vh_t VV \\
(B_{\frac \rho4}(x_{i}), g_{i}) @>\op{GH}>> (B_{\frac \rho4}(x),d_\infty),
\endCD
$$
where $\op{id}_{B_{\frac \rho4}(x_i)}\to h_t$ as $i\to \infty$. It is clear that
$h_t$ is $\Psi(t|n)$-GHA, and by (4.4.1) and (4.4.2), $h_t$ is also a $c(n)$ bi-H\"older homeomorphism. By now we may fix
$t=t_0\in (0,\epsilon(n,\rho,\alpha_0)^{2})$.

For $x\in X$, by (2.5.4) and (2.5.5) we apply the proof of Lemma 2.5 for the case
that $|\op{sec}_{M_i}|\le C(n,\alpha_0,t_0)$ to obtain $\ell_x>0$, and for $i$ large,
a fiber bundle map, $f_i: V_i\to B_{\ell_x}(x)$, with respect to
$g_i(t_0)$ and $V_i$ is a disk bundle over
an infra-nilmanifold $F_i$. Clearly, $h_{t_0}\circ f_i: V_i\to B_{\ell_x}(x)$
gives the desired fiber bundle.
\qed\enddemo

\vskip4mm

\head 5. Proof of Theorem C
\endhead

\vskip4mm

\demo{Proof of Theorem C}

We first show that for any $x\in M$, $\op{vol}(B_\rho(x))\ge v(n,\rho,d)>0$. Arguing by contradiction,
assuming a sequence of collapsing compact $n$-manifolds ($\op{vol}(B_\rho(x_i)\to 0$, equivalently, $\op{vol}(M_i)\to 0$),
$M_i@>\op{GH}>>X$, satisfying
$$\op{Ric}_{M_i}\ge -(n-1),\, d(p_i,q_i)=\op{diam}(M_i)\le d,\,  e_{p_i,q_i}(x_i)\to 0,\, \forall\, x_i\in M_i \text{ is $(\rho,\delta(n))$-LRRP}.$$

Following the almost identical argument in the proof of Theorem 1.1 for $n\ge 4$, one sees that $\pi_1(M_i)\cong \pi_1(V(p_i)\cup_\partial V(q_i))$ is not finite, a contradiction.

To see that $M$ is homeomorphic to a sphere, by a standard compactness argument it
suffices to consider a sequence $M_i@>\op{GH}>>X$ satisfying the above conditions,
and show that $X$ is a homeomorphic $n$-sphere, and for $i$ large, $M_i$ is diffeomorphic to
$X$.

Because $\op{vol}(B_\rho(x_i))\ge v(n,\rho,d)>0$, following the proof of Corollary 1.3 we obtain a constant $c_1(n,\rho,d)>0$ such that any $x_i\in M_i$ is $(c_1(n,\rho,d),\delta(n))$-RP. By (1.4.2), $X$ is a manifold, and for $i$ large $M_i$ is
diffeomorphic to $X$. Because $e_{p,q}(X)\equiv 0$, $X$ is homeomorphic to the gluing of two small
balls at $p, q\in X$, ($d(p,q)=\op{diam}(X)$) that are homeomorphic to Euclidean balls, thus
$X$ is homeomorphic to a sphere.

For $n=3$, $M_i@>\op{GH}>> X$, with $e_{p,q}(X)\equiv 0$. If $\dim(X)=3$, then $X$ is homeomorphic to a sphere, thus $M_i$ is homeomorphic to a sphere. Otherwise, $X$ is isometric to a closed interval, or a suspension surface. By the proof of Theorem 1.1 for $n=3$ (via local Ricci flows in \cite{HW1}), we conclude that $M_i$ is homeomorphic to a gluing of two copies of $D^2\times S^1$ (which may collapse to a half closed and half open interval), thus $M_i$ is homeomorphic to a spherical space form.
\qed\enddemo

\vskip4mm

\noindent {\bf Acknowledgement}: The authors would like to thank Shicheng Xu for helpful comments on Corollaries 0.6 and 0.7, and thanks for some comments from a referee that helps in improving some exposition in this paper.

\enddocument